\documentclass[12pt,leqno]{article}
\usepackage{amssymb}
\usepackage{amsmath}
\usepackage{cite}
\usepackage{amssymb}
\usepackage{amsfonts}
\usepackage{amscd}
\usepackage[matrix]{xypic}

\newtheorem{theorem}{Theorem}

\newtheorem{definition}[theorem]{Definition}

\newtheorem{lemma}[theorem]{Lemma}

\newtheorem{proposition}[theorem]{Proposition}
\newtheorem{remark}[theorem]{Remark}

\newcommand{\A}{\mathcal{A}}

\newcommand{\p}{\mathcal{P}}
\newcommand{\OO}{\mathcal{O}}
\newcommand{\K}{\Bbb{K}}

\begin{document}

\title{\textbf{A class of nonassociative algebras including flexible and alternative algebras, operads and deformations. }}
\author{Elisabeth REMM \thanks{%
corresponding author: e-mail: E.Remm@uha.fr} - Michel GOZE \thanks{%
M.Goze@uha.fr.} \\
\\
{\small Universit\'{e} de Haute Alsace, F.S.T.}\\
{\small 4, rue des Fr\`{e}res Lumi\`{e}re - 68093 MULHOUSE - France}}
\date{}
\maketitle

\begin{abstract}
There exists two types of nonassociative algebras whose associator satisfies a symmetric relation associated 
with a 1-dimensional invariant 
vector space with respect to the natural action of the symmetric group $\Sigma _3$. The first one corresponds to the 
Lie-admissible algebras and this class has been studied in a previous paper \cite{G.R.nonass}. 
Here we are interested by the second one corresponding to the third power associative algebras. 

\end{abstract}

\bigskip


\bigskip


\section{Introduction}

In \cite{G.R.nonass} we have classified, for binary algebras, relations of nonassociativity which are invariant
with respect to an action of the symmetric group on three elements $\Sigma_3$ on the associator. In particular we have
   investigated
two classes of nonassociative algebras, the first one corresponds to algebras whose associator $A_\mu$ satisfies
\begin{eqnarray}
A_\mu-A_\mu \circ \tau_{12}-A_\mu \circ \tau_{23}-A_\mu \circ \tau_{13}+A_\mu \circ c+A_\mu \circ c^2=0,
\end{eqnarray}
and the second
\begin{eqnarray}
\label{p3Ass}
A_\mu+A_\mu \circ \tau_{12}+A_\mu \circ \tau_{23}+A_\mu \circ \tau_{13}+A_\mu \circ c+A_\mu \circ c^2=0.
\end{eqnarray}
where $\tau_{ij}$ denotes the transposition exchanging $i$ and $j$, $c$ is the $3$-cycle $(1,2,3).$ 

These relations are in correspondence with the only two irreducible one dimensional subspaces of $\K[\Sigma_3]$
with respect to the action of $\Sigma_3$, where $\K[\Sigma_3]$ is
the group algebra of $\Sigma_3$. In \cite{G.R.nonass}, we have study, 
the operadic and deformations aspects of the
first one which is the class of Lie-admissible algebras. We will now investigate 
the second class and in particular nonassociative algebras satisfying
(\ref{p3Ass}) with nonassociative relations in correspondence with the subgroups of $\Sigma_3$.

\noindent {\bf Convention:} We consider algebras over a field $\K$ of characteristic zero.

\section{$G_i$-$p^3$-associative algebras}

\subsection{Definition}

Let $\{ G_i \}_{i=1,\cdots,6}$ be the subgroups of $\Sigma_3.$ To fix notations we define
$$\left\{
\begin{array}{l}
G_1=\left\{ Id \right\}, \\
G_2=<\tau_{12}>, \\
G_3=<\tau_{23}>, \\
G_4=<\tau_{13}>, \\
G_5=<c>, \\
G_6=\Sigma_3.
\end{array}
\right.$$
where $<\sigma>$ the cyclic group subgroup generated by $\sigma.$
To each  subgroup $G_i$ we associate the vector $v_{G_i}$ of $\K[\Sigma_3]:$
$$v_{G_i}=\sum_{\sigma\in G_i}\sigma.$$

\begin{lemma}
The one dimensional subspace $\K\{v_{\Sigma_3} \}$ of $\K[\Sigma_3]$ generated by
$$v_{G_6}=v_{\Sigma_3}=\sum_{\sigma \in \Sigma_3} \sigma$$
 is an irreducible invariant subspace of $\K[\Sigma_3].$
\end{lemma}
\begin{definition}
A $G_i$-$p^3$-associative algebra is a $\K$-algebra $(\A,\mu)$ whose associator
$$A_\mu=\mu \circ (\mu \otimes Id - Id \otimes \mu)$$
 satisfies
 $$A_\mu \circ \Phi^{\A}_{v_{G_6}}=0,$$
 where $\Phi^{\A}_{v_{G_i}}:\A^{\otimes^3} \rightarrow \A^{\otimes^3} $ is the linear map
 $$\Phi^{\A}_{v_{G_i}}(x_1\otimes x_2\otimes x_3)=\sum_{\sigma \in G_i}(x_{\sigma^{-1}(1)}\otimes x_{\sigma^{-1}(2)} \otimes x_{\sigma^{-1}(3)} ).$$
\end{definition}
Let $\OO(v_{G_i})$ be the orbit of $v_{G_i}$ with respect to the natural action
$$
\begin{array}{rcl}
\Sigma_3 \times \K[\Sigma_3] &  \rightarrow & \K[\Sigma_3] \\
(\sigma , \sum_i a_i\sigma_i) & \mapsto & \sum_i a_i \sigma^{-1} \circ \sigma_i.
\end{array}
$$
Then putting $F_{v_{G_i}}
=\K(\OO(v_{G_i}))$ we have
$$\left\{
\begin{array}{l}
\dim F_{v_{G_1}}=6, \\
\dim F_{v_{G_2}}=\dim F_{v_{G_3}}=\dim F_{v_{G_4}}=3, \\
\dim F_{v_{G_5}}=2, \\
\dim F_{v_{G_6}}=1. \\
\end{array}
\right.
$$
\begin{proposition}
Every $G_i$-$p^3$-associative algebra is third power associative.
\end{proposition}
Recall that a third power associative algebra is an algebra $(\A,\mu)$ whose associator satisfies $A_\mu(x,x,x)=0.$ 
Linearizing this relation, we obtain
$$A_\mu \circ \Phi^{\A}_{v_{\Sigma_3}}=0.$$ Since each of the invariant spaces $F_{G_i}$ contains the vector $v_{\Sigma_3}$, we deduce the proposition.

\begin{remark}
An important class of third power associative algebras is the class of power associative algebras, that is, algebras such that any element generates an associative
subalgebra.
\end{remark}

\subsection{What are $G_i$-$p^3$-associative algebras?}

\begin{enumerate}
\item If $i=1$, since $v_{G_1}=Id$, the class of $G_1$-$p^3$-associative algebras is the full class of associative algebras.
\item If $i=2$, the associator of a $G_2$-$p^3$-associative algebra $\A$ satisfies
$$A_\mu(x_1,x_2,x_3)+A_\mu(x_2,x_1,x_3)=0$$
and this is equivalent to
$$A_\mu(x,x,y),$$
for all $x,y \in \A.$
\item If $i=3$, the associator of a $G_3$-$p^3$-associative algebra $\A$ satisfies
$$A_\mu(x_1,x_2,x_3)+A_\mu(x_1,x_3,x_2)=0,$$
that is,
$$A_\mu(x,y,y),$$
for all $x,y \in \A$.

Sometimes $G_2$-$p^3$-associative algebras  are called left-alternative algebras, $G_3$-$p^3$-associative algebras  
are  right-alternative algebras.
An alternative algebra is  an algebra which satisfies the $G_2$ and $G_3$-$p^3$-associativity.
\item If $i=4$, we have $A_\mu(x,y,x),$ for all $x,y \in \A,$ and the class of $G_3$-$p^3$-associative algebras is the class of flexible algebras.
\item If $i=5$, the class of $G_5$-$p^3$-associative algebras corresponds to $G_5$-associative algebras 
(\cite{G.R.Lieadm}).
\item If $i=6$, the  associator of a $G_6$-$p^3$-associative algebra satisfies
$$
\begin{array}{l}
A_\mu(x_1,x_2,x_3)+A_\mu(x_2,x_1,x_3)+A_\mu(x_3,x_2,x_1)\\+A_\mu(x_1,x_3,x_2)+A_\mu(x_2,x_3,x_1)
A_\mu(x_3,x_1,x_2)=0.
\end{array}$$
If we consider the symmetric product $x\star y=\mu(x,y)+\mu(y,x)$ and the skew-symmetric product $[x,y]=\mu(x,y)-\mu(y,x)$, then the $G_6$-$p^3$-associative identity
is equivalent to $$[x\star y,z]+ [y\star z,x]+[z\star x,y]=0.$$
\begin{definition}
A $([\, , \, ] , \star)$-admissible-algebra is a $\K$-vector space $\A$ provided  with two multiplication:
\begin{enumerate}
\item a symmetric multiplication $\star$,
\item a skew-symmetric multiplication $[,]$
satisfying the identity
$$[x\star y,z]+ [y\star z,x]+[z\star x,y]=0$$
for any $x,y \in \A$.
\end{enumerate}
\end{definition}

Then a $G_6$-$p^3$-associative algebra can be defined as $([\, , \, ] , \star)$-admissible algebra.
\end{enumerate}

\section{The operads $G_i$-$p^3\A ss $ and their dual}

For each $i\in \left\{ 1, \cdots , 6 \right\}$, the operad for $G_i$-$p^3$-associative 
algebras will be denoted by $G_i$-$p^3\A ss$.
The operads $\left\{ G_i \mbox{ \rm -} p^3\A ss \right\}_{i=1,\cdots,6}$ are binary quadratic operads, that is, operads of the form $\p =\Gamma(E)/(R)$, where $\Gamma(E)$ denotes the free operad generated by a
$\Sigma_2$-module $E$ placed in arity 2 and $(R)$ is the operadic ideal generated by a $\Sigma_3$-invariant subspace $R$
of $\Gamma(E)(3).$ Then the dual operad  $\p^!$ is the quadratic operad $\p^! := \Gamma (E^\vee)/(R^\perp)$,
where $R^\perp \subset \Gamma (E^\vee)(3)$ is the annihilator of $R \subset \Gamma(E)(3)$ in the pairing
\begin{eqnarray}
\label{pairing}
\left\{
\begin{array}{l}
<(x_i \cdot x_j)\cdot x_k,(x_{i'} \cdot x_{j'})\cdot x_{k'}>=0, \ {\rm if} \ \{i,j,k\}\neq \{i',j',k'\}, \\
<(x_i \cdot x_j)\cdot x_k,(x_i \cdot x_j)\cdot x_k>=(-1)^{\varepsilon(\sigma)},  \\
\qquad \qquad \qquad   {\rm with} \ \sigma =
\left(
\begin{array}{lll}
i&j&k \\
i'&j'&k'
\end{array}
\right)
\ {\rm if} \ \{i,j,k\}= \{ i',j',k'\}, \\
<x_i \cdot (x_j\cdot x_k),x_{i'} \cdot (x_{j'}\cdot x_{k'})>=0, \ {\rm if} \ \{ i,j,k\} \neq \{ i',j',k'\}, \\
<x_i \cdot (x_j\cdot x_k),x_i \cdot (x_j\cdot x_k)>=-(-1)^{\varepsilon(\sigma)}, \\
\qquad \qquad \qquad  {\rm with} \ \sigma =
\left(
\begin{array}{lll}
i&j&k \\
i'&j'&k'
\end{array}
\right)
\ {\rm if} \ (i,j,k)\neq (i',j',k'), \\
<(x_i \cdot x_j)\cdot x_k,x_{i'} \cdot (x_{j'}\cdot x_{k'})>=0, 
\end{array}
\right.
\end{eqnarray}
and $(R^\perp)$ is the
operadic ideal generated by $R^\perp.$
For the general notions of binary quadratic operads see \cite{G.K,M.S.S}. 
Recall that a quadratic operad $\p$ is Koszul if the free $\p$-algebra based on a $\K$-vector space $V$ is Koszul, 
for any vector space $V$. This property is conserved by duality and can be studied using generating functions of $\p$ 
and of $\p^!$ (see \cite{G.K} or \cite{M.R.nonKoszul}) Before  studying the Kozsulness of the operads $G_i$-$p^3\A ss$, 
we will compute the homology of an associative algebra which will be useful to look if $G_i$-$p^3\A ss$ are Kozsul or not.

\medskip

\noindent Let $A_2$ the two dimensional associative algebra given in a basis $\{e_1,e_2\}$ by $e_1e_1=e_2,e_1e_2=e_2e_1=e_2e_2=0.$ Recall that the Hochschild homology of an associative algebra is given by the complex $(C_n(\mathcal{A},\mathcal{A}),d_n)$ where
$C_n(\mathcal{A},\mathcal{A})=\mathcal{A}\otimes \mathcal{A}^{\otimes n }$ and the differentials
$d_n:C_n(\mathcal{A},\mathcal{A}) \rightarrow C_{n-1}(\mathcal{A},\mathcal{A}) $
are given by
\[
\begin{array}{ll}
\medskip
d_n(a_0,a_1, \cdots , a_n)= & (a_0a_1,a_2,\cdots a_n)+
\sum_{i=1}^{n-1}(-1)^i(a_0,a_1, \cdots , a_ia_{i+1}, \cdots ,a_n) \\
 & +(-1)^n(a_na_0, a_1, \cdots a_{n-1}).
\end{array}
\]
Concerning the algebra $A_2$, we have
$$d_1(e_i,e_j)=e_ie_j-e_je_i=0,$$
for any $i,j=1,2.$ Similarly we have
$$
\left\{
\begin{array}{l}
d_2(e_1,e_1,e_1)=2(e_2,e_1)-(e_1,e_2),\\
d_2(e_1,e_1,e_2)=d_2(e_1,e_2,e_1)=-d_2(e_2,e_1,e_1)=(e_2,e_2)
\end{array}
\right.
$$
and $0$ in all the other cases. Then $\dim Imd_2=2$ and $\dim Kerd_1=4.$ Then $H_1(A_2,A_2)$ is isomorphic to $A_2$. We have also
$$
\left\{
\begin{array}{ll}
d_3(e_1,e_1,e_1,e_1)&=-(e_1,e_2,e_1)+(e_1,e_1,e_2),\\
d_3(e_1,e_1,e_1,e_2)&=(e_2,e_1,e_2)-(e_1,e_2,e_2),\\
d_3(e_1,e_1,e_2,e_1)&=-d_2(e_2,e_1,e_1,e_1)=(e_2,e_2,e_1)-(e_2,e_1,e_2),\\
d_3(e_1,e_2,e_1,e_1)&=(e_1,e_2,e_2)-(e_2,e_2,e_1),\\
d_3(e_1,e_1,e_2,e_2)&=-d_3(e_1,e_2,e_2,e_1)=-d_3(e_2,e_1,e_1,e_2)=d_3(e_2,e_2,e_1,e_1)\\
&=(e_2,e_2,e_2)
\end{array}
\right.
$$
and $d_3=0$ in all the other cases. Then $\dim Imd_3=4$ and $\dim Kerd_2=6$. 
Thus $H_2(A_2,A_2)$ is non trivial and $A_2$ is not a Koszul algebra.

\medskip

\noindent Now we will study all the operads $G_i$-$p^3\A ss.$

\medskip

\subsection{The operad $(G_1 \mbox{\rm -}p^3\A ss )$ } 
Since $G_1$-$p^3 \A ss=\A ss$, where $\A ss$ denotes the operad for associative algebras, and since the operad
$\A ss$ is selfdual, we have
$$(G_1 \mbox{\rm -}p^3\A ss)^!=\A ss^! =G_1 \mbox{\rm -}p^3\A ss. $$
We also have
$$\widetilde{G_1 \mbox{\rm -}p^3\A ss}=\widetilde{\A ss}=\A ss,$$
where $\widetilde{\p}$ is the maximal current operad of $\p$ defined in \cite{G.R.current}.

\medskip

\subsection{The operad $(G_2 \mbox{\rm -}p^3\A ss )$ }
The operad $G_2$-$p^3 \A ss$ is the operad for left-alternative algebras. It is the quadratic operad
$\p =\Gamma (E)/(R)$, where  the $\Sigma_3$-invariant subspace $R$
of $\Gamma(E)(3)$ is generated by the vectors
$$(x_1\cdot x_2)\cdot x_3-x_1\cdot (x_2\cdot x_3)+(x_2\cdot x_1)\cdot x_3-x_2\cdot (x_1\cdot x_3).$$
The annihilator $R^\perp$ of $R$ with respect to the pairing (\ref{pairing}) is generated by the vectors
\begin{eqnarray}
\left\{
\begin{array}{l}
(x_1 \cdot x_2) \cdot x_3-x_1 \cdot(  x_2 \cdot x_3), \\
(x_1 \cdot x_2) \cdot x_3+(x_2 \cdot  x_1) \cdot x_3.
\end{array}
\right.
\end{eqnarray}
We deduce from direct calculations that $\dim R^\perp=9$ and
\begin{proposition}
The $(G_2 \mbox{\rm -}p^3\A ss )^!$-algebras are associative algebras satisfying
$$abc=-bac.$$
\end{proposition}
Recall that $(G_2 \A ss )^!$-algebras are associative algebras satisfying
$$abc=bac.$$ and this operad is classically denoted $\p erm.$

\begin{theorem}\cite{Dzu}
The operad $(G_2 \mbox{\rm -}p^3\A ss )^!$ is not Koszul.
\end{theorem}

\noindent {\it Proof}.
It is easy to describe $(G_2 \mbox{\rm -}p^3\A ss )^!(n)$ for any $n$. In fact
$(G_2 \mbox{\rm -}p^3\A ss )^!(4)$ correspond to associative elements satisfying
$$x_1x_2x_3x_4=-x_2x_1x_3x_4=-x_2(x_1x_3)x_4=x_1x_3x_2x_4=-x_1x_2x_3x_4$$
and $(G_2 \mbox{\rm -}p^3\A ss )^!(4)=\left\{ 0 \right\}.$ Let $\p$ be $(G_2 \mbox{\rm -}p^3\A ss ).$
The generating function of $\p^!=(G_2 \mbox{\rm -}p^3\A ss )^!$ is
$$g_{\p^!}(x)=\sum_{a \geq 1} \frac{1}{a!}\dim (G_2 \mbox{\rm -}p^3\A ss )^!(a)x^a=
x+x^2+\frac{x^3}{2}.$$
But the generating function of $\p=(G_2 \mbox{\rm -}p^3\A ss )$ is
$$g_{\p}(x)=x+x^2+\frac{3}{2}x^3+\frac{5}{2}x^4+O(x^5)$$
and if $(G_2 \mbox{\rm -}p^3\A ss )$ is Koszul, then the generating functions should be related by the functional equation
$$g_{\p}(-g_{\p ^!}(-x))=x$$
and it is not the case so both  $(G_2 \mbox{\rm -}p^3\A ss )$ and $(G_2 \mbox{\rm -}p^3\A ss )^!$ are not
Koszul.

\begin{remark}
Since the redaction of this work, Dzumadildaev and Zusmanovich have shown and published
this result. For this reason we refer this theorem to these authors.
\end{remark}

By definition, a quadratic operad $\p$ is Koszul if any free $\p$-algebra on a vector space $V$
is a Koszul algebra. Let us describe the free algebra $\mathcal{F}_{(G_2 \mbox{\rm -}p^3\A ss )^!}(V)$
when $\dim V=1$ and $2$.

A $(G_2 \mbox{\rm -}p^3\A ss )^!$-algebra $\mathcal{A}$ is an associative algebra satisfying
$$xyz=-yxz,$$
for any $x,y,z \in \mathcal{A}.$ This implies $xyzt=0$ for any $x,y,z \in \mathcal{A}.$
In particular we have
\[
\left\{
\begin{array}{l}
x^3=0, \\
x^2y=0,
\end{array}
 \right.
\]
for any $x,y \in \mathcal{A}.$
If $\dim V=1$, $\mathcal{F}_{(G_2 \mbox{\rm -}p^3\A ss )^!}(V)$ is of dimension $2$ and given by
\[
\left\{
\begin{array}{l}
e_1e_1=e_2, \\
e_1e_2=e_2e_1=e_2e_2=0.
\end{array}
 \right.
\]
In fact if $V=\K\left\{ e_1 \right\}$ thus in $\mathcal{F}_{(G_2 \mbox{\rm -}p^3\A ss )^!}(V)$
we have $e_1^3=0.$ We deduce that $\mathcal{F}_{(G_2 \mbox{\rm -}p^3\A ss )^!}(V)=A_2$ and
$\mathcal{F}_{(G_2 \mbox{\rm -}p^3\A ss )^!}(V)$ is not Koszul. It is easy to generalize this construction. 
If $\dim V=n$, then
$\dim \mathcal{F}_{(G_2 \mbox{\rm -}p^3\A ss )^!}(V)=\frac{n(n^2+n+2)}{2}$ and if 
$\left\{ e_1, \cdots , e_n \right\}$ is a basis of $V$ then $\left\{ e_i, e_i^2, e_ie_j, e_le_me_p \right\},$ 
for $i,j=1, \cdots ,n$ and $l,m,p=1,\cdots ,n$
with $m>l$, is a basis of $\mathcal{F}_{(G_2 \mbox{\rm -}p^3\A ss )^!}(V).$ For example, if $n=2$, 
the basis of $\mathcal{F}_{(G_2 \mbox{\rm -}p^3\A ss )^!}(V)$ is
$$\left\{ v_1=e_1, v_2=e_2, v_3=e_1^2, v_4=e_2^2, v_5=e_1e_2,v_6=e_2e_1,v_7=e_1e_2e_1,
v_8=e_1e_2^2 \right\}$$ and the multiplication table is
\[
\begin{array}{c|c|c|c|c|c|c|c|c}
    & v_1 & v_2 & v_3 & v_4 & v_5 & v_6 & v_7 & v_8 \\
\hline
v_1 & v_3 & v_5 &  0  & v_8 &  0  & v_7 &  0  &  0  \\
\hline
v_2 & v_6 & v_4 &-v_7 &  0  &-v_8 &  0  &  0  &  0  \\
\hline
v_3 &  0  &  0  &  0  &  0  &  0  &  0  &  0  &  0  \\
\hline
v_4 &  0  &  0  &  0  &  0  &  0  &  0  &  0  &  0  \\
\hline
v_5 & v_7 & v_8 &  0  &  0  &  0  &  0  &  0  &  0  \\
\hline
v_6 &-v_7 &-v_8 &  0  &  0  &  0  &  0  &  0  &  0  \\
\hline
v_7 &  0  &  0  &  0  &  0  &  0  &  0  &  0  &  0  \\
\hline
v_8 &  0  &  0  &  0  &  0  &  0  &  0  &  0  &  0  \\
\end{array}
\]
 For this algebra we have
 \[
 \left\{
 \begin{array}{l}
 d_1(v_1 ,v_2)=v_5-v_6, \\
  \frac{1}{2}d_1(v_1, v_6)=v_7=-d_1(v_1, v_5)=d_1(v_2,v_3), \\
  \frac{1}{2}d_1(v_2, v_5)=-v_8=d_1(v_6, v_2)=-d_1(v_1,v_4), \\
 \end{array}
 \right.
 \]
and $Ker \, d_1$ is of $\dim 64$. The space $Im \, d_2$ doesn't contain in particular the vectors
$(v_i,v_i)$ for $i=1,2$ because these vectors $v_i$ are not in the derived subalgebra. 
Since these vectors are in $Ker \, d_1$
we deduce that the second space of homology is not trivial.  \begin{proposition}
The current operad of $G_2 \mbox{\rm -}p^3\A ss$ is
$$\widetilde{G_2 \mbox{\rm -}p^3\A ss}=\p erm.$$
\end{proposition}
This is directly deduced of the definition of the current operad \cite{G.R.current}.
\subsection{The operad $(G_3 \mbox{\rm -}p^3\A ss )$ }

It is defined by the module of relations generated by the vector
$$(x_1x_2)x_3-x_1(x_2x_3)+ (x_1x_3)x_2-x_1(x_3x_2),$$
and $R^{\perp}$ is the linear span of
\[
\left\{
\begin{array}{l}
(x_1x_2)x_3-x_1(x_2x_3), \\
(x_1x_2)x_3+(x_1x_3)x_2.
\end{array}
\right.
\]
\begin{proposition}
A $(G_3 \mbox{\rm -}p^3\A ss )^!$-algebra is an associative algebra $\A$ satisfying
$$abc=-acb,$$
for any $a,b,c \in \A.$
\end{proposition}

Since $(G_3 \mbox{\rm -}p^3\A ss )^!$ is basically isomorphic to $(G_2 \mbox{\rm -}p^3\A ss )^!$
we deduce that $(G_3 \mbox{\rm -}p^3\A ss )$ is not Koszul.

\subsection{The operad $(G_4 \mbox{\rm -}p^3\A ss )$ }

Remark that a $(G_4 \mbox{\rm -}p^3\A ss )$-algebra is generally called flexible algebra.
The relation
$$A_{\mu} (x_1,x_2,x_3)+ A_{\mu} (x_3,x_2,x_1)=0$$
is equivalent to $A_{\mu} (x,y,x)=0$ and this denotes the flexibility of $(\mathcal{A},\mu).$

\begin{proposition}
A $(G_4 \mbox{\rm -}p^3\A ss )^!$-algebra is an associative algebra satisfying
$$abc=-cba.$$
\end{proposition}

This implies that
$$\dim (G_4 \mbox{\rm -}p^3\A ss )^!(3)=3.$$
We have also $x_{\sigma(1)}x_{\sigma(2)}x_{\sigma(3)}x_{\sigma(4)}=(-1)^{\varepsilon (\sigma)}x_1x_2x_3x_4$ for any $\sigma \in \Sigma _4.$ This gives $\dim (G_4 \mbox{\rm -}p^3\A ss )^!(4)=1.$
Similarly
\[
\begin{array}{ll}
x_1x_2(x_3x_4x_5)&= -x_3(x_4x_5x_2)x_1=x_1(x_4x_5(x_2x_3)) = - x_1x_2(x_3x_5x_4)\\
& =(x_1x_2(x_4x_5))x_3=- (x_4x_5)(x_2x_1)x_3 = (x_3x_2x_1)x_4x_5 \\
& =-x_1x_2x_3x_4x_5
\end{array}
\]
(the algebra is associative so we put some parenthesis just to explain how we pass from one expression
to an other). We deduce $(G_4 \mbox{\rm -}p^3\A ss )^!(5)=\left\{ 0 \right\}$ and more generally $(G_4 \mbox{\rm -}p^3\A ss )^!(a)=\left\{ 0 \right\}$ for $a \geq 5.$

The generating function of  $(G_4 \mbox{\rm -}p^3\A ss )^!$ is
$$f(x)=x+x^2+\frac{x^3}{2}+\frac{x^4}{12}.$$
Let $\mathcal{F}_{(G_4 \mbox{\rm -}p^3\A ss )^!}(V)$ be the free $(G_4 \mbox{\rm -}p^3\A ss )^!$-algebra based on the vector space $V.$ In this algebra we have  the relations
\[
\left\{
\begin{array}{l}
a^3=0, \\
aba=0,
\end{array}
\right.
\]
for any $a,b \in V.$ Assume that $\dim V=1.$ If $\left\{ e_1 \right\}$ is a basis of $V$, then
$e_1^3=0$ and $\mathcal{F}_{(G_4 \mbox{\rm -}p^3\A ss )^!}(V)=\mathcal{F}_{(G_2 \mbox{\rm -}p^3\A ss )^!}(V).$
We deduce that $\mathcal{F}_{(G_4 \mbox{\rm -}p^3\A ss )^!}(V)$ is not a Koszul algebra.

\begin{proposition}
The operad for flexible algebra is not Koszul.
\end{proposition}

Let us note that, if $\dim V=2$ and $\left\{ e_1,e_2 \right\}$ is a basis of $V,$
then $\mathcal{F}_{(G_4 \mbox{\rm -}p^3\A ss )^!}(V)$ is generated by
$\left\{ e_1,e_2,e_1^2,e_2^2, e_1e_2,e_2e_1, e_1e_2^2, e_1^2e_2, e_2e_1^2,e_2^2e_1,e_1^2e_2^2,
e_2^2e_1^2\right\}$ and is of dimension $12$.
\begin{proposition}
We have
$$\widetilde{G_4 \mbox{\rm -}p^3\A ss}=(G_4 \mbox{\rm -}\A ss )^!.$$
\end{proposition}
This means that a $\widetilde{G_4 \mbox{\rm -}p^3\A ss}$ is an associative algebra $\A$ satisfying
$$abc=cba,$$
for any $a,b,c\in \A.$
\subsection{The operad $(G_5 \mbox{\rm -}p^3\A ss )$ }

It coincides with $(G_5 \mbox{\rm -}\A ss )$ and this last has been studied in \cite{G.R.Lieadm}.

\subsection{The operad $(G_6 \mbox{\rm -}p^3\A ss )$ }

A $(G_6 \mbox{\rm -}p^3\A ss )$-algebra $(\mathcal{A},\mu)$ satisfies the relation
\[
\begin{array}{c}
A_{\mu }(x_1,x_2,x_3)+A_{\mu }(x_2,x_1,x_3)+A_{\mu }(x_3,x_2,x_1)\\ +A_{\mu }(x_1,x_3,x_2)
+A_{\mu }(x_2,x_3,x_1)+A_{\mu }(x_3,x_1,x_2)=0.
\end{array}
\]
The dual operad $(G_6 \mbox{\rm -}p^3\A ss )^!$ is generated by the relations
\[
\left\{
\begin{array}{l}
(x_1x_2)x_3=x_1(x_2x_3), \\
(x_1x_2)x_3=(-1)^{\varepsilon (\sigma)}(x_{\sigma(1)} x_{\sigma(2)}) x_{\sigma(3)}, \ \rm{for \ all \ } \sigma \in \Sigma_3.
\end{array}
\right.
\]
We deduce
\begin{proposition}
A $(G_6 \mbox{\rm -}p^3\A ss )^!$-algebra is an associative algebra $\A$ which satisfies
$$abc=-bac=-cba=-acb=bca=cab, $$
for any $a,b,c \in \mathcal{A}$. In particular
\[
\left\{
\begin{array}{l}
a^3=0,\\
aba=aab=baa=0.
\end{array}
\right.
\]
 \end{proposition}

\begin{lemma} The operad $(G_6 \mbox{\rm -}p^3\A ss )^!$ satisfies
$(G_6 \mbox{\rm -}p^3\A ss )^!(4)=\left\{0 \right\}.$
\end{lemma}

\noindent {\it Proof.} We have in $(G_6 \mbox{\rm -}p^3\A ss )^!(4)$ that
\[
x_1(x_2x_3)x_4=x_2(x_3x_4)x_1=-x_1(x_3x_4x_2)=x_1x_3x_2x_4=-x_1x_2x_3x_4
\]
so $x_1x_2x_3x_4=0.$ We deduce that the generating function of $(G_6 \mbox{\rm -}p^3\A ss )^!$ is
$$f^!(x)=x+x^2+\frac{x^3}{6}.$$
If this operad is Koszul the generating function of the operad $(G_6 \mbox{\rm -}p^3\A ss )$ should be
of the form
$$f(x)=x+x^2+\frac{11}{6}x^3+\frac{25}{6}x^4+\frac{127}{12}x^5+\cdots$$
But if we look the free algebra generated by $V$ with $\dim V=1$, it satisfies $a^3=0$ and coincides
with $\mathcal{F}_{(G_2 \mbox{\rm -}p^3\A ss )^!}(V).$ Then $(G_6 \mbox{\rm -}p^3\A ss )$ is not Koszul.
\begin{proposition}
We have
$$\widetilde{G_6 \mbox{\rm -}p^3\A ss}=\mathcal{L}ieAdm^!$$ that is the binary quadratic operad whose corresponding algebras
are associative and satisfying
\[
abc=acb=bac.
\]
\end{proposition}

\section{Cohomology and Deformations}

Let $(\A,\mu)$ be a $\K$-algebra defined by quadratic relations. 
It is attached to a quadratic linear operad $\p$. By deformations of $(\A,\mu)$, we mean (\cite{G.R.valued})

\begin{itemize}

\item a $\K^*$ non archimedian extension field of $\K$,
with a valuation $v$ such that, if
$A$ is the ring of valuation and $\mathcal{M}$ the unique ideal of $A$, then the residual field $A/\mathcal{M}$ is isomorphic to $\K$.

\item The $A/\mathcal{M}$ vector space $\widetilde{\A}$ is $\K$-isomorphic to $\A$.

\item For any $a,b \in \A$ we have that
$$\tilde{\mu}(a,b)-\mu(a,b)$$
belongs to the $\mathcal{M}$-module $\widetilde{\A}$ (isomorphic to $\A \otimes \mathcal{M}$)
.

\end{itemize}
The most important example concerns the case where $A$ is $\K[[t]],$ the ring of formal series. In this case $\mathcal{M}=\left\{ \sum_{i \geq 1}a_it^i,a_i \in \K \right\},$ $\K^*=\K((t))$
the field of rational fractions. This case corresponds to the classical Gerstenhaber deformations. 
Since $A$ is a local ring, all the notions of valued deformations coincides (\cite{Fialowsky}).

We know (\cite{Ma}) that there exists always a cohomology which parametrizes  deformations. 
If the operad $\p$ is Koszul, this cohomology is the "standard''-cohomology called the operadic cohomology. 
If the operad $\p$ is not Koszul, the cohomology which governs deformations is based on the minimal model of 
$\p$ and the operadic cohomology and deformations cohomology differ.

In this section we are interested by the case of left-alternative algebras, that is, by the operad
$(G_2 \mbox{\rm -}p^3\A ss )$ and also by the classical alternative algebras.

\subsection{Deformations and cohomology of left-alternative algebras}

A $\K$-left-alternative algebra $(\A, \mu)$ is a $\K$-$(G_2 \mbox{\rm -}p^3\A ss )$-algebra. Then $\mu$ satisfies
$$A_\mu(x_1,x_2,x_3)+A_\mu(x_2,x_1,x_3)=0.$$
A valued deformation can be viewed as a $\K[[t]]$-algebra $(A \otimes \K[[t]],\mu_t)$ whose product
$\mu _t$ is given by
$$\mu _t =\mu +\sum_{i \geq 1} t^i\varphi _i.$$

\medskip

a) {\it The operadic cohomology}

\medskip

\noindent It is the standard cohomology $H^*_{(G_2 \mbox{\rm -}p^3\A ss )}(\A,\A)_{st}$ of the 
$(G_2 \mbox{\rm -}p^3\A ss )$-algebra $(\A,\mu)$. 
It is associated to the cochains complex
$$\mathcal{C}^1_{\p}(\A,\A)_{st} \xrightarrow{\delta_{st} ^1} \mathcal{C}^2_{\p}(\A,\A)_{st} \xrightarrow{\delta_{st} ^2}
\mathcal{C}^3_{\p}(\A,\A)_{st} \xrightarrow{\delta_{st} ^3} \cdots $$
where $\p=(G_2 \mbox{\rm -}p^3\A ss )$ and 
$$\mathcal{C}^p_{\p}(\A,\A)_{st} =Hom(\p^!(p) \otimes_{ \Sigma_p} \A ^{ \otimes p},\A).$$
Since $(G_2 \mbox{\rm -}p^3\A ss )^!(4)=0,$ we deduce that
$$ H^p_{\p}(\A,\A)_{st}=0 \ {\rm for} \ p \geq 4,$$
because the cochains complex is a short sequence
$$\mathcal{C}^1_{\p}(\A,\A)_{st} \xrightarrow{\delta_{st} ^1} \mathcal{C}^2_{\p}(\A,\A)_{st} \xrightarrow{\delta_{st} ^2}
\mathcal{C}^3_{\p}(\A,\A)_{st} \xrightarrow{0} 0.$$
The coboundary operator are given by
\[\left\{
\begin{array}{cll}
\delta^1f(a,b)&=& f(a)b +af(b)-f(ab), \\
\delta^2\varphi(a,b,c)&=& \varphi(ab,c)+\varphi(ba,c)-\varphi(a,bc)-\varphi(b,ac)\\
& & \varphi(a,b)c+\varphi(b,a)c-a\varphi(b,c)-b\varphi(a,c).
\end{array}
\right.
\]

\medskip

b) {\it The deformations cohomology}

\medskip

\noindent The minimal model of $G_2 \mbox{\rm -}p^3\A ss $ is a homology isomorphism
$$(G_2 \mbox{\rm -}p^3\A ss,0) \xrightarrow{\rho} (\Gamma(E),\partial)  $$
of dg-operads such that the image of $\partial$ consists of
decomposable elements of the free operad $\Gamma (E).$ Since $G_2 \mbox{\rm -}p^3\A ss(1)=\K$,
this minimal model exists and it is unique. The deformations cohomology $H^*(\A,\A)_{defo}$
of $\A$ is the cohomology of the complex
$$\mathcal{C}^1_{\p}(\A,\A)_{defo} \xrightarrow{\delta ^1} \mathcal{C}^2_{\p}(\A,\A)_{defo} \xrightarrow{\delta ^2}
\mathcal{C}^3_{\p}(\A,\A)_{defo} \xrightarrow{\delta ^3} \cdots $$
where
\[
\left\{
\begin{array}{l}
\mathcal{C}^1_{\p}(\A,\A)_{defo}=Hom(\A,\A),\\
\mathcal{C}^k_{\p}(\A,\A)_{defo}=Hom(\oplus_{q\geq 2} E_{k-2}(q) \otimes_{ \Sigma_q} \A ^{ \otimes q},\A).
\end{array}
\right.
\]
The Euler characteristics of $E(q)$ can be read off from the inverse of the generating function of the operad 
$G_2 \mbox{\rm -}p^3\A ss $
$$g_{G_2 \mbox{\rm -}p^3\A ss }(t)=t+t^2+\frac{3}{2}t^3
+\frac{5}{2}t^4+\frac{53}{12}t^5$$
which is
$$g(t)=t-t^2+\frac{t^3}{2}
+\frac{13}{3}t^5+O(t^6).$$
We obtain in particular
\[
\begin{array}{l}
\chi(E(4))=0.
\end{array}
\]
Each one of the modules $E(p)$ is a graded module
$(E_*(p))$ and
$$\chi (E(p))=\dim E_0(p)-\dim E_1(p)+\dim E_2(p)+\cdots$$
We deduce
\begin{itemize}
\item[--]
$E(2)$ is generated by two degree $0$ bilinear operation $\mu_2 : V \cdot V \to V$,
\item[--]
$E(3)$ is generated by three degree $1$ trilinear operation $\mu_3 :
    V^{\otimes ^3} \to V$,
\item[--]
$E(4)=0$.
\end{itemize}

Considering the action of $\Sigma_n$ on $E(n)$ we deduce that $E(2)$ is generated by a binary operation of degree 0 whose differential satisfies
$$\partial (\mu_2)=0,$$
 $E(3)$ is generated by a trilinear operation of degree one such that
$$\partial (\mu_3)=\mu_2 \circ_1 \mu_2-\mu_2 \circ_2 \mu_2+\mu_2 \circ_1 (\mu_2 \cdot \tau_{12})-(\mu_2 \circ_2 \mu_2)\cdot \tau_{12}.$$

(we have $(\mu_2 \circ_2 \mu_2)\cdot \tau_{12}(a,b,c)=b(ac)$)

Since $E(4)=0$ we deduce

\begin{proposition}
The cohomology $H^*(\A,\A)_{defo}$ which governs deformations or right-alternative algebras is associated to the complex
\[\mathcal{C}^1_{\p}(\A,\A)_{defo} \xrightarrow{\delta ^1} \mathcal{C}^2_{\p}(\A,\A)_{defo} \xrightarrow{\delta ^2}
\mathcal{C}^3_{\p}(\A,\A)_{defo} \xrightarrow{\delta ^3} \mathcal{C}^4_{\p}(\A,\A)_{defo} \rightarrow  \cdots
\]
with $$
\begin{array}{l}
\mathcal{C}^1_{\p}(\A,\A)_{defo}=Hom(V^{\otimes 1},V),  \\
\mathcal{C}^2_{\p}(\A,\A)_{defo}=Hom(V^{\otimes 2},V),  \\
\mathcal{C}^3_{\p}(\A,\A)_{defo}=Hom(V^{\otimes 3},V),  \\
\mathcal{C}^4_{\p}(\A,\A)_{defo}=Hom(V^{\otimes 5},V)\oplus \cdots \oplus Hom(V^{\otimes 5},V),
\end{array}
$$
In particular any $4$-cochains consists of $5$-linear maps.
\end{proposition}

\subsection{Alternative algebras}

Recall that an alternative algebra is given by the relation
\[
A_\mu(x_1,x_2,x_3)=-A_\mu(x_2,x_1,x_3)=A_\mu(x_2,x_3,x_1).
\]
\begin{theorem}
An algebra $(\A,\mu)$ is alternative if and only if the associator satisfies
\[
A_\mu \circ \Phi_v^{\A}=0,
\]
with $v=2Id+\tau_{12}+\tau_{13}+\tau_{23}+c_1$.
\end{theorem}

\noindent {\it Proof.} The associator satisfies $A_\mu \circ \Phi_{v_1}^{\A}=A_\mu \circ \Phi_{v_2}^{\A}$ with
$v_1=Id +\tau_{12}$ and $v_2=Id +\tau_{23}$. The invariant subspace of $\K[\Sigma_3]$ generated by $v_1$ and $v_2$ is of
dimension $5$ and contains the vector $\sum_{\sigma \in \Sigma_3}\sigma.$ From \cite{G.R.nonass}, the space is generated by the orbit of the vector $v$.

\begin{proposition}
Let $\A lt$ be the operad for alternative algebras. Its dual is the operad for associative algebras satisfying
$$abc-bac-cba-acb+bca+cab=0.$$
\end{proposition}

\noindent {\bf Remark.} The current operad $\widetilde{\A lt}$ is the operad for associative algebras satisfying
$abc=bac=cba=acb=bca$, that is, $3$-commutative algebras so
$$\widetilde{\A lt}=\mathcal{L}ie\A dm^! \, .$$

In \cite{Dzu}, one gives the generating functions of $\p=\A lt$ and $\p^!=\A lt^!$
\[
\begin{array}{l}
g_{\p}(x)=x+\frac{2}{2!}x^2+\frac{7}{3!}x^3+\frac{32}{4!}x^4+\frac{175}{5!}x^5+\frac{180}{6!}x^6+O(x^7),\\
g_{\p^!}(x)=x+\frac{2}{2!}x^2+\frac{5}{3!}x^3+\frac{12}{4!}x^4+\frac{15}{5!}x^5.
\end{array}
\]
and conclude to the non Koszulity of $\A lt$.

The operadic cohomology is the cohomology associated to the complex
\[
(\mathcal{C}^p_{\A lt}(\A,\A)_{st} =(Hom(\A lt^!(p)\otimes_{\Sigma_p} \A^{\otimes p},\A), \delta_{st}).
\]
Since $\A lt^!(p)=0$ for $p\geq 6$ we deduce the short sequence
\[
\mathcal{C}^1_{\A lt}(\A,\A)_{st} \xrightarrow{\delta_{st} ^1} \mathcal{C}^2_{\A lt}(\A,\A)_{st} \rightarrow \cdots
\rightarrow \mathcal{C}^5_{\A lt}(\A,\A)_{st} \rightarrow 0.
\]
But if we compute the formal inverse of the function $-g_{\A lt}(-x)$ we obtain
\[
x+x^2+\frac{5}{6}x^3+\frac{1}{2}x^4+\frac{1}{8}x^5-\frac{11}{72}x^6+O(x^7).
\]
Because of the minus sign it can not be the generating function of the operad $\p ^!=\A lt^!.$ So this implies also that both operad are not Koszul.
But it gives also some information on the deformation cohomology. In fact if $\Gamma (E)$ is the free operad associated to the minimal model, then
$$
\left\{
\begin{array}{l}
\dim \chi (E(2))=-2,\\
\dim \chi (E(3))=-5,\\
\dim \chi (E(4))=-12,\\
\dim \chi (E(5))=-15,\\
\dim \chi (E(6))=+110.\\
\end{array}
\right.
$$
Since
 $\chi (E(6))=\sum _i (-1)^i \dim E_i(6)$, the graded space $E(6)$ is not concentred in degree even. Then the $6$-cochains of the deformation cohomology are
$6$-linear maps of odd degree.

\end{document}